\documentclass[12pt]{amsart}
\usepackage{amsmath, amsthm, amscd, amsfonts, amssymb}

\usepackage{amsmath}
\usepackage{amsthm}
\usepackage{amssymb}
\usepackage{amscd}
\usepackage[mathscr]{euscript}

\newtheorem{Lemma}{Lemma}
\newtheorem{Theorem}[Lemma]{Theorem}
\newtheorem{Proposition}[Lemma]{Proposition}
\newtheorem{Corollary}[Lemma]{Corollary}
\newtheorem{Definition}[Lemma]{Definition}
\newtheorem{Example}[Lemma]{Example}
\newtheorem{Remark}[Lemma]{Remark}

\newcommand{\bthe}{\begin{Theorem}}
	\newcommand{\ethe}{\end{Theorem}}
\newcommand{\bpr}{\begin{Proposition}}
	\newcommand{\epr}{\end{Proposition}}
\newcommand{\bco}{\begin{Corollary}}
	\newcommand{\eco}{\end{Corollary}}
\newcommand{\bre}{\begin{Remark} \em}
	\newcommand{\ere}{\end{Remark}}
\newcommand{\ble}{\begin{Lemma}}
	\newcommand{\ele}{\end{Lemma}}
\newcommand{\bex}{\begin{Example} \em}
	\newcommand{\eex}{\end{Example}}
\newcommand{\bde}{\begin{Definition} \em}
	\newcommand{\ede}{\end{Definition}}
\newcommand{\bpf}{\begin{proof}[\bf{Proof.}]}
	\newcommand{\epf}{\end{proof}}
\newcommand{\bqu}{\begin{Question} \em}
	\newcommand{\equ}{\end{Question}}

\newcommand{\qq}{\mathbb{Q}}
\newcommand{\zz}{\mathbb{Z}}
\newcommand{\nn}{\mathbb{N}}
\newcommand{\ff}{\mathbb{F}}
\newcommand{\End}{{\rm End}}
\newcommand{\Hom}{{\rm Hom}}

\newcommand{\Soc}{{\rm Soc}}
\newcommand{\ann}{{\rm ann}}

\newcommand{\mat}{{\rm M}}
\newcommand{\car}{{\rm char}}

\begin{document}

\title{On Pure subRings of sp-Groups}

\author[A. Amini]{A. Amini}
\address{Afshin Amini, Department of Mathematics, College of Sciences, Shiraz University, Shiraz 71457, Iran}
\email{aamini@shirazu.ac.ir}

\author[B. Amini]{B. Amini}
\address{Babak Amini, Department of Mathematics, College of Sciences, Shiraz University, Shiraz 71457, Iran}
\email{bamini@shirazu.ac.ir}

\author[E. Momtahan]{E. Momtahan}
\address{Ehsan Momtahan, Department of Mathematics,
	Yasouj University, Yasouj, Iran}
\email{e-momtahan@mail.yu.ac.ir}


	
\begin{abstract}
	Let $G$ be a sp-group such that for every prime $p$, $G_p$ is elementary. 
	We show that $\End_{\zz}(G)$ is a sp-group and  every subring $R$ of $\prod \End_{\zz}(G_p)$, containing $\oplus \End_{\zz}(G_p)$ is pure  if and only if $R=\mathbb{M}_T=\{x\in \prod_{p\in \mathbb{P}}\End(G_p) \;|\; \exists k\in \nn \;\mbox{\rm{such that}} \;\; kx \in T \},$ where $T$ is a subring of $\prod_{p\in \mathbb{P}}\End(G_p)$. We observe that $\frac{\mathbb{M}_T}{\oplus_{p\in \mathbb{P}}\End(G_p)}$ is (ring) isomorphic with $T\otimes_{\zz} \qq$.  Moreover, we conclude  that  a significant number of the examples around the topic can be easily obtained and described by choosing an appropriate subring $T$.  
\end{abstract}


\maketitle

 Following  \cite{ktt}, a reduced group $G$ with  infinitely many nonzero $p$-components satisfying the equivalent conditions of the next proposition is called an sp-group. 
\bpr {\rm{\cite[Proposition 1.4]{ktt}}} Let $G$ be a reduced mixed group which has infinitely many nonzero $p$-components. The following properties of the group $G$ are equivalent.
\begin{enumerate}
\item  For every prime $p$, we have the direct decomposition $G = G_p\oplus B_p$ for some group $B_p$ with $pB_p = B_p$.

\item  The embeddings $$\oplus G_p < G \leq \prod G_p,$$
hold and $G$ is pure subgroup of $\prod G_p$.

\item  The embeddings from item (2) hold and $\frac{G}{G_t}$ is a divisible group.
\item  The $p$-component $G_p$ is a direct summand of the group $G$ for every prime $p$ and $\frac{G}{G_t}$ is a
divisible group.
\end{enumerate}
\epr

In this article we deal mainly with sp-groups with elementary $p$-components, that is, for every prime $p$, $G_p$ is a direct sum of cyclic groups of order $p$. 
By group we always mean an Abelian group and $\prod$ and $\oplus$, denote $\prod_{p\in \mathbb{P}}\End(G_p)$ and $\oplus_{p\in \mathbb{P}}\End(G_p)$ respectively, unless otherwise stated. Since for every prime $p$,  $G_p$ is  an elementary $p$-group, $\End_{\zz}(G_p)$ is (ring) isomorphic with $\End_{\zz_p}(V)$, where $V$ is a (right) vector space over $\zz_p$. Hence $\prod$ is a regular right self-injective ring. When for every $p$, $G_p$ is finite, or equivalently  $\End(G_p)=\mat_{n_p}(\zz_p)$, $\prod$ is right and left self-injective. Also $G_t$ denotes the torsion subgroup of $G$ and  $G_t=\oplus_{p\in \mathbb{P}}G_p$. A subring $R$ of $\prod$ is said to be pure if its additive group is a pure subgroup of $\prod$. We say that a subgroup $H$ of a group $G$ is pure if for every $n\in \nn$, $nH=H\cap nG$. By a regular ring we always mean regular in the sense of von Neumann (for a systematic study of regular rings see \cite{go}). A ring $R$ is said to be right $\aleph_0$-self-injective ($p$-injective)
if for every (module) homorphism $f\in \Hom_R(I,R)$, there
exists $\bar{f}\in \Hom_R(R,R)$, such that $\bar{f}|I = f$, where $I$ is any countably
generated (principal) right ideal of $R$. The reader is referred to \cite{la} and \cite{fu}, for undefined terms and notations.

In the sequel we need a lemma, which has partially generalized some parts of the aforementioned proposition. 

\ble \label{sp} Let $G$ be an Abelian group, such that $\oplus G_p \leq G \leq \prod G_p$ and  for every $p$, $G_p$ is elementary. Then $G$ is a pure subgroup of $\prod G_p$  if and only if  for every $p$, $G=G_p\oplus pG$.

\ele
\bpf
 Let $G$ be pure, and $p$ is a prime number. It is clear that $\prod G_p=G_p\oplus \prod_{q\neq p} G_q=G_p\oplus p\prod G_p$. Now by modular law and the fact that $G$ is pure, we have $G=G_p\oplus (G\cap p\prod G_p)=G_p\oplus pG$. Conversely, suppose that for every prime $p$, $G=G_p\oplus pG$. This immediately implies that $pG=p^2G$ and hence $pG=p^nG$ for every $n\in \nn$. For two distinct prime numbers $p$ and $q$, we have $G=G_p\oplus G_q\oplus pqG$. This can be easily generalized to any finite number of prime numbers. Now suppose that $n=p_1^{\alpha_1}\cdots p_k^{\alpha_k}$ is a natural number, then $G\cap n\prod G_p=nG\oplus (G_{p_1}\oplus \cdots \oplus G_{p_n}\cap n\prod G_p)=nG$.
\epf

\bpr \label{sp} Let $G$ be a sp-group, such that  for every prime $p$, $G_p$ is elementary. Then 
$\End_{\zz}(G)$ is a sp-group.	


\epr
\bpf 
 Since for every prime $p$, $G=G_p\oplus pG$ and $pG=p^2G$, we conclude that  $\Hom_{\zz}(G_p,pG)$ and $\Hom(pG,G_p)$ are both zero. Therefore $\End_{\zz}(G)=\End_{\zz}(G_p)\oplus \End_{\zz}(pG)$. Thus, $\oplus \End_{\zz}(G_p)$ is contained in $\End_{\zz}(G)$. Now consider the short exact sequence $0\to G_t \to G \to \frac{G}{G_t}\to 0$ and apply the functor $\Hom(-,G)$, we get  $$\End_{\zz}(G)\subseteq \Hom_{\zz}(G_t,G)=\End_{\zz}(G_t)=\prod \End_{\zz}(G_p)$$ (recall that $\Hom_{\zz}(\frac{G}{G_t}, G)=0$ due to $\frac{G}{G_t}$ being divisible and $G$ being reduced). Now we show that $\End_{\zz}(G)$ is a pure subring of $\prod \End_{\zz}(G)_p$. Since $\End(G)_p=\End(G_p)$, we conclude that $\End(G)\leq \prod \End(G)_p$. Using Lemma \ref{sp}, it is enough to show that $\End_{\zz}(G)=\End_{\zz}(G)_p\oplus p\End_{\zz}(G)$ for every prime $p$. According to the equality $\End(G)=\End(G_p)\oplus \End(pG)$, we have to show that $\End(pG)=p\End(G)$. We know that $p\End(G)=p\End(pG)$. Hence we show that $\End(pG)=p\End(pG)$. Suppose that $f\in \End(pG)$ and $x\in pG$. We see that $f(x)=py$ for some $y\in pG$, due to this fact that $pG=p^2G$. Now define $h:pG\longrightarrow pG$ with $h(x)=y$. The map $h$ is well-defined, because if $f(x)=py=pz$ for some $y,z\in pG$, then $p(y-z)=0$, i.e., $y-z=0$ (remind that $pG\cap G_p=0$). That is $f=ph$, for $h\in \End(pG)$, this proves that $\End(pG)=p\End(pG)$. Therefore, $\End(G)$ is a pure subring of $\prod\End(G)_p$.

\epf
Note that under the assumptions of Proposition \ref{sp}, $\End(G)$ is a pure subring of $\prod$. In the following we characterize the format of all pure subrings of $\prod$.

\bde Let $T$ be a subset of $\prod$. Then by $\mathbb{M}_T$, we mean 
$$\{x\in \prod\;|\; \exists k\in \nn \;\mbox{\rm{such that}} \;\; kx \in  \oplus + T \}.$$
If, in addition, $T$ satisfies in this property that when $x\in T$, then $kx\in T$ for every $k\in \zz$ (e.g., $T$ is a subgroup of $\prod$), then $\mathbb{M}_T$ can be equivalently defined as  $$\mathbb{M}_T=\{x\in \prod\;|\; \exists k\in \nn \;\mbox{\rm{such that}} \;\; kx \in  T \}.$$
\ede

\bthe \label{shape}
Let $T$ be a subring of $\prod$, 
then the following hold: 

\begin{enumerate}
\item $\mathbb{M}_T$ is a pure subring of $\prod$,
containing $\oplus + T.$

\item $\frac{\mathbb{M}_T}{\oplus}$ is (ring) isomorphic with $T\otimes \qq.$

\item $\mathbb{M}_T$ is a regular ring if and only if $T\otimes \qq$ is a regular ring.

{\rm{If in addition, for every $p$, $G_p$ is finite, the following facts hold either}}:
\item $\mathbb{M}_T$ is right and left non-singular.

\item $\prod$ is a right (a left) maximal quotient ring of $\mathbb{M}_T$.

\item $\mathbb{M}_T$ is never right or left $\aleph_0$-self-injective ring unless $\mathbb{M}_T=\prod$.

\item $|\mathbb{M}_T|=|T|$.



\end{enumerate}

\ethe
\bpf
(1): Let $x,y\in \mathbb{M}_T$, we know that there exist $k,l\in \nn$ such that  $kx=s$ and $ly=t$, where $s,t\in T$. We observe that $kl(xy)=(kx)(ly)=st\in  T$. On the other hand $kl(x+y)=l(kx)+k(ly)=ls+kt$ which belong to $T$ as well. This shows that $\mathbb{M}_T$ is a ring. We claim that $\mathbb{M}_T$  is a pure subgroup of $\prod$, i.e., for $n\in \nn$, we have $\mathbb{M}_T\cap n\prod=n\mathbb{M}_T$. To show this, suppose that $x\in \mathbb{M}_T\cap n\prod$, so $x=ny$ for some $y\in \prod$ and there exists $k\in \nn$ such that $kx\in T$, so $kny\in T$, hence $y\in \mathbb{M}_T$, so $x\in n\mathbb{M}_T$.\

(2) Define $\phi:\mathbb{M}_T\longrightarrow T\otimes \qq$ with $\phi(x)=t\otimes \frac{1}{k}$, where $k\in \nn$, $t\in T$ and $kx=t\in T$. First of all, we show that $\phi$ is well-defined. Suppose that for $x\in \mathbb{M}_T$, there are $k,k'\in \nn$ such 
$kx=t$ and $k'x=t'$. Note that $k'kx=k't=kt'$, which implies that $k't\otimes \frac{1}{kk'}=kt'\otimes \frac{1}{kk'}$.  Hence $t\otimes \frac{1}{k}=t'\otimes \frac{1}{k'}$, i.e., $\phi$ is well-defined.
Now we show that $\phi$ is a ring homomorphism. To show that $\phi(x+y)=\phi(x)+\phi(y)$, suppose that for $x,y\in \mathbb{M}_T$, there are $k,l\in \nn$ such that $kx=s$ and $ly=t$. Now consider 
$lkx=ls$ and $kly=kt$,we have 
$kl(x+y)=ls+kt$, this implies that 
$\phi(x+y)=ls+kt\otimes\frac{1}{kl}$, but $ls\otimes\frac{1}{kl}+kt\otimes \frac{1}{kl}=s\otimes \frac{1}{k}+t\otimes \frac{1}{l}=\phi(x)+\phi(y)$. Now, since $klxy=(kx)(ly)=st$, we have $\phi(xy)=st\otimes\frac{1}{kl}=(s\otimes\frac{1}{k})(t\otimes \frac{1}{l})=\phi(x)\phi(y)$. Furthermore, we show that  $\phi$ is onto. Without loss of generality we may suppose that  $s\otimes \frac{1}{k}\in T\otimes \qq$, now the equation $k(x+\oplus)=s+\oplus$ is solvable, due to $\prod/\oplus$ being divisible. Hence, there exists $x\in \prod$ such that $kx=s+a$, where $a\in \oplus$. There exists $l\in \nn$ such that $la=0$,  so $lkx=ls$, therefore $x\in \mathbb{M}_T$ and $\phi(x)=ls\otimes \frac{1}{lk}=s\otimes\frac{1}{k}$. Now we show that $\ker \phi=\oplus$. To observe this, recall first  that $T\otimes \qq\cong S^{-1}T$, where $S=\zz\setminus\{0\}$. Now  $t\otimes \frac{1}{k}=0$ if and only if there exists $n\in \zz\setminus\{0\}$ such that $nt=0$, i.e., $t\in  \oplus$, due to $G_p$ is elementary for every prime $p$, i.e.,  $\ker \phi=\oplus$. \

(3) Since $\oplus$ is regular, we deduce that $\mathbb{M}_T$ is regular if and only if $T\otimes \qq$ is a regular ring.\

(4) Since $\ann_{\prod}(\oplus)\cap\mathbb{M}_T=\ann_{\mathbb{M}_T}(\oplus)$ and $\ann_{\prod}(\oplus)=(0)$, we have $\ann_{\mathbb{M}_T}(\oplus)=(0)$. But $\Soc(\mathbb{M}_T)=\oplus$, hence $\mathcal{Z}(\mathbb{M}_T)\subseteq \ann(\Soc(\mathbb{M}_T))=(0)$ (see \cite[Proposition 7.13]{la}). Thus $\mathbb{M}_T$ is right and left non-singular.\ 

(5)  Since 
$\oplus$ is an essential ideal of $\prod$, and by (4), $\mathbb{M}_T$ is non-singular, we are thorough by \cite[Proposition 13.39,(2)]{la}.\

(6) Let $\mathbb{M}_T$ be a proper subring of $\prod$. By (4),  $\mathbb{M}_T$ is non-singular and $\Soc(\mathbb{M}_T)=\oplus$ is countably generated. If $\mathbb{M}_T$ were $\aleph_0$-self-injective, it were self-injective, due to \cite[Corollary 4.4]{mo}, but this is not the case because the maximal quotient ring of $\mathbb{M}_T$ is $\prod$.\

(7) In as much as  $|\frac{\mathbb{M}_T}{\oplus}|=|T\otimes \qq|=|T|$ and $|\oplus|=\aleph_0$, we have $|\mathbb{M}_T|=\aleph_0 |T|=|T|$.

\epf

\bre \label{direct} (i) Let $T$ be a subring of $\prod$ and $I$ be a right ideal of $T$; then $\mathbb{M}_I$ is a right ideal of $\mathbb{M}_T$, containing $\oplus$ and $\frac{\mathbb{M}_T}{\mathbb{M}_I}\cong \frac{T\otimes \qq}{I\otimes \qq}\cong \frac{T}{I}\otimes \qq$.\\
(ii) In general, it is not the case that every pure subring of $\prod$, contains $\oplus$. For example, consider $R=\prod T_2(\zz_p)$ which is a pure subring of $\prod \mat_2(\zz_p)$, where by $T_2(\zz_p)$ we mean the subring of all upper triangular matrices of $\mat_2(\zz_p)$. However, the next result shows that those pure subrings of $\prod$, which contain $\oplus$, have the $\mathbb{M}_T$ shape.

\ere

\bpr Let $R$ be a subring of $\prod$ containing $\oplus$, then $R$ is a pure subring if and only if there exists a subring $T$ of $\prod$ such that $R=\mathbb{M}_T$. 

\epr
\bpf
($\Leftarrow$) it has already proved in Theorem \ref{shape}, part 1.\\
($\Rightarrow$): Let $R$ be a pure subring of $\prod $, we show that $\mathbb{M}_R= R $. It is clear $R$ is a subring of $\mathbb{M}_R$. Now let $x\in \mathbb{M}_R$, then there exists $k\in \nn$ such that $kx\in \oplus + R = R$. Hence $kx\in (R\cap k\prod)$, but $R$ is pure, hence $R\cap k\prod=kR$, which implies that $kx\in kR$, i.e., $kx=kr$, for $r\in R$. Thus $k(x-r)=0$, or equivalently  $x-r\in \oplus \subset R$, and hence $x-r, r\in R$, so $x\in R$.

\epf

\section*{The Curious Case of $\prod_{p\in \mathbb{P}}\zz_p$}
In the last few decades, one of the main sources of  examples and counter-examples  of sp-groups was $\prod_{p\in \mathbb{P}}\zz_p$. This gives us permission to examine this particular but important case closely. One of the goals of this section is to take a fresh look  at some old examples in this field. We will also construct other examples while scrutinizing the old ones. 




The first close attentions to regular subrings of $\prod_{p\in \mathbb{P}}\zz_p$, containig $\oplus_{p\in \mathbb{P}}\zz_p$, very likely, goes back to \cite{fh}, where the authors proved their embedding theorem, i.e., every regular ring is embedded in a regular ring with identity. To do this, they constructed a regular ring $\mathbb{M}$ with $1$ as follows. Let  $\varepsilon_p$, be the identity
	of $\zz_p$ and $\varepsilon=(\cdots, \varepsilon_p, \cdots)$ be the identity of $\prod_{p\in \mathbb{P}}\zz_p$. Now consider the quotient $\frac{\prod_{p\in \mathbb{P}}\zz_p}{\oplus_{p\in \mathbb{P}}\zz_p}$, which is a torsion-free divisible ring in which
	the pure subgroup generated by the coset of $\varepsilon$ is a ring $\frac{\mathbb{M}}{\oplus_{p\in \mathbb{P}}\zz_p}$, isomorphic
	to $\qq$. In this way, they obtained a subring $\mathbb{M}$ of $\prod_{p\in \mathbb{P}}\zz_p$ which contains $\varepsilon$ as its identity and contains
	every $\zz_p$. This $\mathbb{M}$ is regular: it contains the regular ring $\oplus_{p\in \mathbb{P}}\zz_p$, as an ideal
	modulo which $\mathbb{M}$ is regular and they proved every regular ring is a unital
	$\mathbb{M}$-algebra. For more information on $\mathbb{M}$, see \cite{fh}, page 286. In Example \ref{F-H}, we  observe that $\mathbb{M}$ is but only a first step towards a large class of pure subrings.  Before that we need a simple but useful lemma.


\ble \label{exist} Let $\alpha$ be a cardinal number which is  less or equal than $\omega$. Then $\zz^{\alpha}$ is contained in $\prod_{p\in \mathbb{P}} \zz_p$ as a subring.
\ele

\bpf We first provide a proof for the case $\alpha=\omega$. 
We may partition $\mathbb{P}$ into $\aleph_0$ parts, each of which has infinite elements. Suppose that $\mathbb{P}=\bigcup_{i=1}^{\infty}A_i$, and $A_i\cap A_j=\emptyset$, when $i\neq j$ and $|A_i|=\aleph_0$ for every $i\in \nn$. Now we have 

$$\prod_{p\in \mathbb{P}}\zz_p \cong \prod_{p\in A_1}\zz_p \times \prod_{p\in A_2}\zz_p\times \cdots$$
But $\zz$ is contained in $\prod_{p\in A_i}\zz_p$ for every $i\in \nn$, which means that $\zz^{\omega}\subset \prod_{p\in \mathbb{P}}\zz_p$. The proof for $\alpha < \omega$ is almost the same.
\epf

The following Proposition has been first appeared in \cite[Lemma 2]{ra}. We give a proof for the sake of completeness. Note that, in contrast to Remark \ref{direct}, part (ii), it is not difficult to show that a pure subring of $\prod_{p\in \mathbb{P}}\zz_p$ always contains $\oplus_{p\in \mathbb{P}}\zz_p$.

\bpr \label{end}
Let $\mathbb{M}$ be a subring of $\prod_{p\in \mathbb{P}}\zz_p$, containing $\oplus_{p\in \mathbb{P}}\zz_p$ and for every $p\in \mathbb{P}$, $\mathbb{M}=\mathbb{M}_p\oplus p\mathbb{M}$ (i.e., $\mathbb{M}$ is a pure subring of  $\prod_{p\in \mathbb{P}}\zz_p$), then $\End_{\zz}(\mathbb{M})=\mathbb{M}$.
\epr
\bpf
We show that the map $\phi:\mathbb{M}\longrightarrow \End_{\zz}(\mathbb{M})$ with $\phi(a)=\lambda_a$ is an isomorphism, where $\lambda_a:\mathbb{M}\longrightarrow \mathbb{M}$ is defined by $\lambda_a(x)=ax$. It is clear that $\phi$ is a monomorphism. We show that $\phi$ is onto. Let $f\in \End_{\zz}(\mathbb{M})$, we know that 
$$\End_{\zz}(\mathbb{M})\subseteq \Hom(\mathbb{M},\prod_{p\in \mathbb{P}}\zz_p)=\prod_{p\in \mathbb{P}} \Hom(\mathbb{M}, \mathbb{M}_p).$$
Recall that here $\mathbb{M}_p=\zz_p$. Now we define $f_p:\mathbb{M}\longrightarrow \zz_p$ by $f_p:=\pi_ p\circ f$ ($\pi_p:\prod_{p\in \mathbb{P}}\zz_p\longrightarrow \zz_p)$. So for every $x \in \mathbb{M}$, $f(x)=((f_p(x)))_{p\in \mathbb{P}}$. By hypothesis $\mathbb{M}=\zz_p \oplus p\mathbb{M}$, so 
$\Hom(\mathbb{M},\zz_p)=\Hom(\zz_p\oplus p\mathbb{M}, \zz_p)\cong
\Hom(\zz_p, \zz_p)\oplus \Hom(p\mathbb{M}, \zz_p)\cong \Hom(\zz_p, \zz_p)\cong \zz_p.$ Therefore, for every $f_p\in \Hom(\mathbb{M}, \zz_p)$, there exists $a_p\in \zz_p$ such that $f_p((x_p))=a_px_p$. We observe that $f(x)=f((x_p)_{p\in \mathbb{P}})=(f_p(x))_{p\in \mathbb{P}}=(a_p)(x_p)=\lambda_a(x)$, where $a=(a_p)_{p\in \mathbb{P}}$. Since $f(1)=\lambda_a(1)=a$, we conclude that  $a\in \mathbb{M}$ and $\phi$ is onto, i.e., $\End_{\zz}(\mathbb{M})\cong \mathbb{M}$.

\epf

\bre All pure subrings $\mathbb{M}_T$, which 
have been introduced in this section are subject to Proposition \ref{end}, and therefore $\End(\mathbb{M}_T)=\mathbb{M}_T$.

\ere


\bex \label{F-H} Going back to Fusch-Halperin's example, 
we easily observe that $\mathbb{M}=\{x\in \prod_{p\in \mathbb{P}}\zz_p\;|\; \exists k\in \nn \;\mbox{such that} \;\; kx \;\; \mbox{is eventually constant} \}.$ Looking again at $\mathbb{M}$, we may redefine it as $\mathbb{M}_{\zz}=\{x\in \prod_{p\in \mathbb{P}}\zz_p\;|\; \exists k\in \nn \;\mbox{such that} \;\; kx \in  \zz \}.$ In general, let $n\in \nn$ be given, we may define $\mathbb{M}_{\zz^n}=\{x\in \prod_{p\in \mathbb{P}}\zz_p\;|\; \exists k\in \nn \;\mbox{such that} \;\; kx \in  \zz^n \}.$ Furthermore, we define $\mathbb{M}_{\zz^{\omega}}$ as $\{x\in \prod_{p\in \mathbb{P}}\zz_p\;|\; \exists k\in \nn \;\mbox{such that} \;\; kx \in  \zz^{\omega} \}.$ Just  for simplicity's sake, we use notations $\mathbb{M}_n$ and $\mathbb{M}_{\omega}$ instead of $\mathbb{M}_{\zz^n}$ and $\mathbb{M}_{\zz^{\omega}}$ respectively. 
A  reason for regularity of $\mathbb{M}_n$ is that $\frac{\mathbb{M}_n}{\oplus_{p\in \mathbb{P}}\zz_p}\cong \zz^n\otimes \qq \cong \qq^n$. However,  $\mathbb{M}_{\omega}$ is not regular. By Theorem \ref{shape}, 
$\frac{\mathbb{M}_{\omega}}{\oplus_{p\in \mathbb{P}}\zz_p}	\cong \zz^{\omega} \otimes \qq$. 
The right hand side ring is isomorphic to a subalgebra of $\qq^{\omega}$
consisting of those sequences of rational numbers whose denominators are bounded (with respect to suitable (not any) representations as fractions).	That is  $\zz^{\omega}\otimes \qq \cong A=\{\;x\in \qq^\omega\;|\; x=(\;\frac{b_1}{k}, \;\frac{b_2}{k}, \;\frac{b_3}{k}, \cdots) \;\mbox{for} \;k\in \nn\}$, which is not regular. To see this, It is enough to consider the element $a=(\; \frac{1}{2}, \; \frac{3}{2}, \; \frac{5}{2},\; \frac{7}{2}, \cdots)$. It is clear that there is no element $l\in A$ such that $ala=a$. Moreover, it is  worth mentioning that $\mathbb{M}_{\omega}$ is not even a p-injective ring, otherwise it were regular due to $\mathbb{M}_{\omega}$ being reduced.
\eex

\bre \label{ziad}
(I) Let $n\geq 2$, we can partition $\mathbb{P}$ into $n$ (infinite) subsets, in exactly $2^{\aleph_0}$ ways. According to each of these partitions, we have an $\mathbb{M}_n$. The sum of all these $\mathbb{M}_n$'s for a fixed $n$, i.e., $\sum \mathbb{M}_n$, is a regular Baer subring of 
$\prod_{p\in \mathbb{P}}\zz_p$, which contains all idempotents of it (see Example \ref{rang}).\\
(II) In general we cannot say that $\mathbb{M}_n$ is a subring of $\mathbb{M}_k$, where $k\geq n$, but there are ways of partitioning $\mathbb{P}$, in which $\mathbb{M}_1 \subset \mathbb{M}_2 \subset \cdots$. Also there are uncountably many chains like this.

\ere


\bex \label{rang}
 Rangaswamy in \cite{ra} - in a remark that appeared on page 357 of the article - has found a regular Baer subring $B$ of $\prod_{p\in \mathbb{P}}\zz_p$, which is not self-injective. The ring $B$ is the subring generated by $\mathbb{M}_1$ and all idempotents of $\prod_{p\in \mathbb{P}}\zz_p$. Following Theorem \ref{shape}, put $T=\sum \zz e$, where $e$ runs over all idempotents of $\prod_{p\in \mathbb{P}}\zz_p$. Then $B=\mathbb{M}_T$.  In the sequel, we reveal  the relation between $B$ and $\mathbb{M}_n$'s. As we have already mentioned in Remark \ref{ziad}, there is only one $\mathbb{M}_1$ (as there is only one $\zz$ in $\prod_{p\in \mathbb{P}}\zz_p$),  but for every $n\geq 2$, we have uncountably many $\mathbb{M}_n$ (as there are uncountably many $\zz^n$ in $\prod_{p\in \mathbb{P}}\zz_p$). Let $\mathcal{M}_n$ denote the set of all $\mathbb{M}_n$'s. It can be shown that $\sum_{\mathbb{M}_n\in \mathcal{M}_n} \mathbb{M}_n$ is equal to the $B$. Since $\oplus_{p\in \mathbb{P}} \zz_p\oplus \zz^n \subset \mathbb{M}_n$,  we may write $1=\epsilon=\epsilon_1+\cdots+\epsilon_n$, where by $\epsilon_i$ we mean $(0,\cdots,1,\cdots,0)$, where $1$ is in the ith component. So $\mathbb{M}_n=\mathbb{M}_1\epsilon_1+
\cdots+\mathbb{M}_1\epsilon_n$. On the other hand, every idempotent in $\prod_{p\in \mathbb{P}}\zz_p$ belongs to some $\mathbb{M}_n$, for every $n\geq 2$. It is also worth mentioning that $B=\sum_{\mathbb{M}_n\in \mathcal{M}_n} \mathbb{M}_n$ is continuous as well, because it contains all the idempotents of its maximal quotient ring, i.e., $\prod_{p\in \mathbb{P}}\zz_p$ (see \cite[Theorem 13.13]{go}) . On the other hand, we are allowed to write $$\sum_{\mathbb{M}_n\in \mathcal{M}_n} \mathbb{M}_n=\{x\in \prod_{p\in \mathbb{P}}\zz_p\;|\; \exists k\in \nn \;\mbox{such that} \;\; kx \in \sum \zz^n \},$$ where the right hand sum is taken over all (uncountably many) $\zz^n$ in $\prod_{p\in \mathbb{P}}\zz_p$. 
\eex



\bex
We may consider the increasing sequence of regular rings 
$\mathbb{M}_1 < \mathbb{M}_2 < \cdots ,$ then $\varinjlim \mathbb{M}_n$ is a regular subring of $\prod_{p\in \mathbb{P}}\zz_p$ (direct limits of regular rings are regular). The ring is obviously different from $\mathbb{M}_n$'s and $\mathbb{M}_{\omega}$. In fact

$$\varinjlim \mathbb{M}_n=\{x\in \prod_{p\in \mathbb{P}}\zz_p\;|\; \exists k\in \nn \;\mbox{such that} \;\; kx \in \varinjlim\zz^n \}.$$

Since $\varinjlim \zz^n \subset \zz^{\omega}$, we conclude that $\varinjlim \mathbb{M}_n$ is a proper subring of $\mathbb{M}_{\omega}$.
\eex
The following example is another instance of a regular pure subring which is not of finite torsion free rank.
 \bex
Let $T=\zz^{(\omega)}+1.\zz$ be a subring of $\prod_{p\in \mathbb{P}}\zz_p$, then $\frac{\mathbb{M}_T}{\oplus\zz_p}\cong T\otimes \qq \cong \qq^{(\omega)}+ 1\qq$ which is a (regular) subring of $\qq^{\omega}$. Since $\varinjlim\zz^n=T$, we deduce that $\varinjlim \mathbb{M}_n=\mathbb{M}_T$. 
\eex

The next example has been considered in \cite[Example 4.3]{gw} as a pp-subring  (i.e., pricipal ideals are projective) of $\prod_{p\in \mathbb{P}} \zz_p$ which is not regular. In the following, we put the example in the format of $\mathbb{M}_T$, where $T$ is an appropriate subring of $\prod_{p\in \mathbb{P}}\zz_p$. It is worth mentioning that a ``Cantor-like'' proof can be used to show that $\prod_{p\in \mathbb{P}}\zz_p$, has $2^{\aleph_0}$ transcendental elements (non-integral elements) over $\zz$. Hence the number of algebraic (integral) elements over $\zz$ is  $\aleph_0$. But to see a concrete non-integral - in fact transcendental- element in $\prod_{p\in \mathbb{P}}\zz_p$, take the element $d=(1,2,2,2,\cdots)\in \prod_{p\in \mathbb{P}}\zz_p$. Then $\frac{1}{d}=d^{-1}=(1,2^{-1},2^{-1},2^{-1},\cdots)$ exists in $\prod_{p\in \mathbb{P}}\zz_p$. This $\frac{1}{d}$ is actually equal to $(1,2,3,4,6,\cdots ) \in \prod_{p\in \mathbb{P}}\zz_p$. We show that $d^{-1}$ is not integral (and algebraic) over $\zz$. Let $d^{-n}+\sum_{i=1}^n b_id^{-n+i}=0$, which implies that $1+\sum_{i=1}^nb_i d^i=0$. That is $1+\sum_{k=1}^n2^kb_k=0\in \zz_p$ for $p\geq 3$. For those prime numbers $p$ which are enough large it is not possible. The same method shows that $d^{-1}$ is transcendental.
\bex 
We may redefine the ring which has been introduced in \cite[Example 4.3]{gw} as  follows: put $T=\zz[d]$, and define
$$\mathbb{M}_{T}=\{x\in \prod_{p\in \mathbb{P}}\zz_p\;|\; \exists k\in \nn \;\mbox{such that} \;\; kx \in  \zz[d] \},$$ 
where $d$ has already been defined in the above (which is transcendental over $\qq$) (see \cite[Example 4.3]{gw}). According to Theorem \ref{shape}, we have: $$\frac{\mathbb{M}_T}{\oplus\zz_p}\cong T\otimes \qq  \cong \zz[x]\otimes \qq\cong \qq[x],$$ which is clearly not regular.
\eex

\end{document}